\documentclass[11pt]{article}
\usepackage{latexsym}
\newtheorem{thm}{Theorem}
\newtheorem{cor}[thm]{Corollary}

\newtheorem{rem}[thm]{Remark}
\newtheorem{exm}[thm]{Example}
\newtheorem{ass}{Assumption}

\def\proof{\noindent{\bf Proof. }}
\def\endproof{\hfill$\Box$\bigskip}
\def\IR{ I\!\!R}

\newcommand{\ind}{\mathbf{1}}

\newcommand{\abs}[1]{\left\vert#1\right\vert}

\newcommand{\A}{\mathcal{A}}

\begin{document}

\title{Crossing Probabilities for Diffusion Processes
with Piecewise Continuous Boundaries}

\author{Liqun Wang\\
University of Manitoba, Department of Statistics\\
Winnipeg, Manitoba, Canada R3T 2N2\\
E-mail: liqun\_wang@umanitoba.ca
\and
Klaus P\"otzelberger\\
University of Economics and Business Administration Vienna\\
Department of Statistics, Augasse 2-6, A-1090 Vienna, Austria\\
E-mail: klaus.poetzelberger@wu-wien.ac.at}

\maketitle

\begin{abstract}
We propose an approach to compute the boundary crossing
probabilities for a class of diffusion processes which can be
expressed as piecewise monotone (not necessarily one-to-one)
functionals of a standard Brownian motion. This class includes many
interesting processes in real applications, e.g.,
Ornstein-Uhlenbeck, growth processes and geometric Brownian motion
with time dependent drift. This method applies to both one-sided and
two-sided general nonlinear boundaries, which may be discontinuous.
Using this approach explicit formulas for boundary crossing
probabilities for certain nonlinear boundaries are obtained, which
are useful in evaluation and comparison of various computational
algorithms. Moreover, numerical computation can be easily done by
Monte Carlo integration and the approximation errors for general
boundaries are automatically calculated. Some numerical examples are
presented.
\end{abstract}

\paragraph{Keywords:}
Boundary crossing probabilities; Brownian motion; diffusion process;
first hitting time; first passage time; Wiener process.

\paragraph{AMS 2000 Subject Classification:}
Primary 60J60, 60J70, Secondary 60J25, 60G40, 65C05


\section{Introduction}

Let $X = \{X_t, t\geq 0\}$ be a diffusion process defined on a
probability space $(\Omega, \A, P)$ and has state space either the
real space $\IR$ or a subinterval of it. In this paper, we are
concerned with the following boundary crossing probability (BCP)
\begin{equation}\label{bcp-general}
P_X\left(a,b,T\right) = P\left(a(t)<X_t<b(t), \forall t\in[0,
T]\right),
\end{equation}
where $T>0$ is fixed, boundaries $a(t)$ and $b(t)$ are real
functions satisfying $a(t) < b(t)$ for all $0<t\le T$ and $a(0) <
x_0 < b(0)$.

The problem of boundary crossing probabilities, or the first passage
time (FPT) distributions, has drawn tremendous amount of attention
in many scientific disciplines. To mention a few, for example, it
arises in biology (Ricciardi et al, 1999), economics (Kr\"amer,
Ploberger and Alt 1988), engineering reliability (Ebrahimi 2005),
epidemiology (Martin-L\"of 1998, Tuckwell and Wan 2000 and Startsev
2001), quantitative finance (Garrido 1989, Roberts and Shortland
1997, Lin 1998, Novikov et al, 2003, Borovkov and Novikov 2005),
computational genetics (Dupuis and Siegmund 2000), seismology
(Michael 2005), and statistics (Doob 1949, Anderson 1960, Durbin
1971, Sen 1981, Siegmund 1986, Bischoff et al 2003, and Zeileis
2004).

Despite its importance and wide applications, explicit analytic
solutions to boundary crossing problems do not exist, except for
very few instances. Traditionally, the mainstream of the research of
nonlinear boundary problems is based on Kolmogorov partial
differential equations for the transition density function, and
focuses on approximate solutions of certain integral or differential
equations for the first-passage time densities. For example, to
obtain the approximate solutions for Brownian motion crossing
continuously differentiable boundaries, the tangent approximation
and other image methods have been used by Strassen (1967), Daniels
(1969, 1996), Ferebee (1982) and Lerche (1986), whereas a series
expansion method has been used by Durbin (1971, 1992), Ferebee
(1983), Ricciardi et al. (1984), Giorno et al. (1989), and Sacerdote
and Tomassetti (1996) to deal with more general diffusion processes.
Recently, Monte Carlo path-simulation methods have been proposed to
numerically compute the FPT density for a general diffusion process
crossing a one-sided constant boundary (Kloeden and Platen 1992,
Giraudo and Sacerdote 1999, and Giraudo, Sacerdote and Zucca 2001).

In principle, all these FPT methods are designed for the
continuously differentiable boundaries only, and they mostly deal
with one boundary problems. However, two-sided and discontinuous
boundary problems arise in many real applications, e.g., in problems
of pricing barrier options in quantitative finance. Furthermore, the
actual numerical computations in FPT methods are either intractable
or give approximate solutions for which the approximation accuracies
are difficult to assess (Daniels 1996, Sacerdote and Tomassetti
1996). The more recent path-simulation methods require heavy
computation and the algorithms are typically very complicated
(Giraudo, Sacerdote and Zucca 2001).

In contrast to the traditional FPT methods, Wang and P\"otzelberger
(1997) proposed an alternative approach, which focuses on the
boundary crossing probabilities (BCP) directly. Using this approach,
Wang and P\"otzelberger (1997) derived an explicit formula for the
BCP for Brownian motion crossing a piecewise linear boundary, and
then used this formula to obtain approximations of the BCP for more
general nonlinear boundaries. The numerical computation can be
easily done by Monte Carlo integration, and the approximation errors
are automatically computed. This approach also allows one to control
the approximation error before the actual computation begins. This
method has been extended to two-boundary problems by Novikov,
Frishling and Kordzakhia (1999), who use recurrent numerical
integration with Gaussian quadrature, and by P\"otzelberger and Wang
(2001), who use Monte Carlo integration. Moreover, it has been
applied in finance (Novikov, Frishling and Kordzakhia 2003, Borovkov
and Novikov 2005), change-point problems (Zeileis 2004), and
engineering reliability (Ebrahimi 2005).

In this paper, we extend the method of Wang and P\"otzelberger
(1997) and P\"otzelberger and Wang (2001) further to a class of
diffusion processes, which can be expressed as piecewise monotone
(not necessarily one-to-one) functionals of a Brownian motion. This
class contains many interesting processes arising in real
applications, e.g., Ornstein-Uhlenbeck processes, growth processes
and the geometric Brownian motion with time-dependent drift.
Moreover, this approach allows us to derive explicit formulas for
BCP for certain nonlinear boundaries, which are useful in evaluation
and comparison of various computational methods and algorithms for
boundary crossing problems.

From methodology point of view, the framework used in this paper is
different from that in the previous works in the literature. Whereas
the previous works are based on the Kolmogorov partial differential
equations for the transition density functions, our approach is
based on the stochastic differential equation (SDE) for the
processes. The SDE approach becomes more and more popular in modern
times, especially in mathematical finance since the publication of
the stimulating works of Black and Scholes (1973) and Merton (1973).
See, e.g., Garrido (1989), Roberts and Shortland (1997), Tuckwell
and Wan (2000), Kou and Wang (2003), Novikov, Frishling and
Kordzakhia (2003), and Ebrahimi (2005).

The paper is organized as follows. Section 2 extends the BCP
formulas of Wang and P\"otzelberger (1997) and P\"otzelberger and
Wang (2001) for Brownian motion to piecewise continuous boundaries.
Section 3 derives the BCP formulas for more general diffusion
processes. Sections 4, 5 and 6 treat the Ornstein-Uhlenbeck
processes, growth processes and geometric Brownian motion
respectively. Finally, section 7 gives numerical examples and
section 8 contains conclusions and discussion.


\section{Brownian motion and piecewise continuous boundaries}

Let $W = \{W_t, t\ge 0\}$ be a standard Brownian motion (Wiener
process) with $EW_t = 0$ and $EW_tW_s = \min(t,s)$. We consider the
BCP $P_W(a,b) = P\left(a(t)< W_t <b(t), \forall t\in[0, T]\right)$.
Wang and P\"{o}tzelberger (1997), Novikov, Frishling and Kordzakhia
(1999), and P\"{o}tzelberger and Wang (2001) have derived explicit
formulas for $P_W(a,b,T)$ for continuous boundaries $a$ and $b$. It
is straightforward to generalize their results to piecewise
continuous boundaries. Throughout the paper we denote by
$(t_i)_{i=1}^n$, $0<t_1<\cdots <t_{n-1}<t_n=T$, a partition of
interval $[0, T]$ of size $n\geq 1$. Further, let $t_0=0$, $\Delta
t_i=t_i-t_{i-1}$, $\alpha_i = a(t_i)$, $\beta_i = b(t_i)$ and
$\delta_i = \beta_i - \alpha_i$. Then completely analog to Wang and
P\"{o}tzelberger (1997) and P\"{o}tzelberger and Wang (2001), we can
establish the following results.
\begin{thm}\label{thm-linear}
Given any partition $(t_i)_{i=1}^n$ of interval $[0,T]$, suppose the
boundaries $a$ and $b$ are linear functions in each subinterval
$(t_{i-1},t_i)$ and are such that, at any $t_i$, $\lim_{t\to
t_i-}a(t) \geq \lim_{t\to t_i+}a(t)$ and $\lim_{t\to t_i-}b(t) \leq
\lim_{t\to t_i+}b(t)$ hold. Then
\begin{equation}\label{bcp-bm} P_W(a,b,T) = Eg(W_{t_1}, W_{t_2},
..., W_{t_n}),
\end{equation}
where function $g(x)$, $x=(x_1,x_2,...,x_n)'$, is defined as
follows.

(1) For one-sided BCP $P_W(-\infty,b,T)$,
\begin{equation}\label{bcp-one-sided}
g(x) = \prod_{i=1}^n\ind_{(-\infty,\beta_i)}(x_i)\left\{1 -
\exp\left[ - \frac{2}{\Delta t_i}(\beta_{i-1}-x_{i-1})(\beta_i-x_i)
\right]\right\},
\end{equation}
where $\ind(\cdot)$ is the indicator function.

(2) For two-sided BCP $P_W(a,b,T)$,
\begin{equation}\label{bcp-two-sided}
g(x) = \prod_{i=1}^n\ind_{(\alpha_i,\beta_i)}(x_i)\left[1-
\sum_{j=1}^\infty h_{ij}(x_{i-1},x_i) \right],
\end{equation}
where
\begin{eqnarray*}
h_{ij}(x_{i-1},x_i) &=& \exp\left\{-\frac{2}{\Delta
t_i}\left[j\delta_{i-1} + (\alpha_{i-1}-x_{i-1})\right]
\left[j\delta_i + (\alpha_i-x_i)\right]\right\} \\
&-& \exp\left\{-\frac{2j}{\Delta t_i} \left[j\delta_{i-1}\delta_i
+ \delta_{i-1}(\alpha_i-x_i)- \delta_i(\alpha_{i-1}-x_{i-1})\right]\right\} \\
&+& \exp\left\{-\frac{2}{\Delta t_i} \left[j\delta_{i-1} -
(\beta_{i-1}-x_{i-1})\right]\left[j\delta_i -(\beta_i-x_i)\right]\right\} \\
&-& \exp\left\{-\frac{2j}{\Delta t_i} \left[j\delta_{i-1}\delta_i
- \delta_{i-1}(\beta_i-x_i) + \delta_i(\beta_{i-1}-x_{i-1})
\right]\right\}.
\end{eqnarray*}
\end{thm}

\paragraph{}
Note that function $h_{ij}(x_{i-1},x_i)$ in (\ref{bcp-two-sided})
consists of exponential functions that decrease very rapidly with
the increase of $j$. In all of our numerical examples, fairly
accurate approximations are achieved by using six terms only.
Moreover, a simplified version of $(\ref{bcp-two-sided})$ is given
in P\"{o}tzelberger and Wang (2001).

The results of Theorem \ref{thm-linear} can be used to approximate
the BCP $P_W(a,b)$ for general nonlinear boundaries $a$ and $b$,
provided they can be sufficiently well approximated by some
piecewise linear functions. More precisely, if the sequence of
piecewise linear functions $a_n(t)\to a(t)$ and $b_n(t)\to b(t)$
uniformly on $[0,T]$, then it follows from the continuity property
of probability measure that
\begin{equation}\label{bcp-approximate}
\lim_{n\to\infty} P_W(a_n,b_n,T) = P_W(a,b,T).
\end{equation}
The accuracy of approximation (\ref{bcp-approximate}) depends on the
partition size $n$. In general, larger $n$ will give more accurate
approximation. For twice continuously differentiable boundaries,
Novikov, Frishling and Kordzakhia (1999) showed that the
approximation errors converge to zero at rate $\sqrt{\log n/n^3}$,
if equally-spaced piecewise linear approximation is used. Later,
P\"otzelberger and Wang (2001) obtained an approximation rate of
$O(1/n^2)$ by using an "optimal" partition of $[0,T]$. More
recently, Borovkov and Novikov (2005) obtained the same rate by
using equally-spaced partition.

In practice, one simple and straightforward way to assess the
approximation accuracy is to calculate the lower and upper bounds of
the BCP using two piecewise linear functions approaching $a$ and $b$
from each side respectively. This method will be demonstrated
through numerical examples in section 7. See Wang and
P\"{o}tzelberger (1997) and P\"{o}tzelberger and Wang (2001) for
more details and examples.

\section{General diffusion processes}

In this and the subsequent sections, we generalize the results of
the previous section to more general diffusion processes satisfying
the following stochastic differential equation (SDE)
\begin{equation}\label{sde-general}
dX_t =  \mu(t,X_t)dt + \sigma(t,X_t)dW_t,\quad X_0 = x_0,
\end{equation}
where the drift $\mu(t,x): [0,\infty)\times\IR\to\IR$ and diffusion
coefficient $\sigma(t,x): [0,\infty)\times\IR\to\IR_+$ are real,
deterministic functions and $\{W_t, t\geq 0\}$ is the standard
Brownian motion (BM) process. In order for the SDE
(\ref{sde-general}) to admit unique solution, throughout this paper
we assume generally that $\mu(t,x)$ and $\sigma(t,x)$ are measurable
functions and satisfy the following Lipschitz and growth conditions:
for any $t\geq 0$,
\begin{equation}\label{Lipschitz}
\abs{\mu(t,x)-\mu(t,y)} + \abs{\sigma(t,x)-\sigma(t,y)} \leq
L\abs{x-y}, \forall x,y\in\IR,
\end{equation}
and for any $x\in\IR$,
\begin{equation}\label{Growth}
\abs{\mu(t,x)} + \abs{\sigma(t,x)} \leq K(1+\abs{x}), \forall
t\geq 0,
\end{equation}
where $L>0$ and $K>0$ are constants. However, since the solution
problem is not a major concern of this paper, in the subsequent
derivations we will not deal with conditions (\ref{Lipschitz}) and
(\ref{Growth}) explicitly. Rather we refer the reader to, e.g.,
Karatzas and Shreve (1991, p289) or Lamberton and Lapeyre (1996,
p50). See also Yamada and Watanabe (1971). In addition, throughout
the paper we make the following assumption, which is necessary for
the application of It\^o's formula.

\begin{ass}\label{ass1}
Let $\{X_t, t\in [0,T]\}$ denote the solution of the SDE
(\ref{sde-general}) with state space an open or closed finite or
infinite interval $\mathcal{S}$. Assume that $\sigma(t,x)>0$ for all
$t\in [0,T]$ and $x\in\mathcal{S}$. Further, assume that there
exists a function $f(t,x)\in C^{1,2}([0,T]\times\mathcal{S})$ (once
continuously differentiable with respect to $t$ and twice with
respect to $x$), such that 1. for $x\in\mathcal{S}$
\[
\frac{\partial f(t,x)}{\partial x}=\frac{1}{\sigma(t,x)}
\]
and 2. It\^o's Lemma may be applied to $Y_t=f(t,X_t)$, i.e.
\begin{equation}\label{Itosform}
dY_t=\frac{\partial f(t,X_t)}{\partial t} dt+\frac{\partial
f(t,X_t)}{\partial x}dX_t +\frac{1}{2}\frac{\partial^2
f(t,X_t)}{\partial x^2}\sigma(t,X_t)^2 dt.
\end{equation}
\end{ass}

Now we consider the BCP (\ref{bcp-general}) for any diffusion
process $X=\{X_t,t\geq 0\}$ satisfying (\ref{sde-general}) and
Assumption \ref{ass1}. The basic idea is that
$P_X\left(a,b,T\right)$ can be computed through the BCP for the
standard BM, as long as $X$ can be expressed as a "nice"
transformation of it. More precisely, given any boundaries $a(t)$,
$b(t)$ and $T>0$, let us denote
\[
A(a,b,T) = \{x(t)\in C([0,\infty)): a(t)<x(t)<b(t), \forall
t\in[0,T]\}.
\]
Then it is easy to see that, if there exists a standard BM $\tilde
W = \{\tilde W_t,t\geq 0\}$ and a measurable functional $f$, such
that $X=f(\tilde W)$, then
\begin{equation}\label{bcp-transform}
P(X\in A(a,b,T)) = P(\tilde W\in f^{-1}A(a,b,T)).
\end{equation}
Consequently, the BCP for $X$ can be calculated through the BCP for
$\tilde W$, as long as there exist boundaries $c(t)$, $d(t)$ and
$S>0$, such that $f^{-1}A(a,b,T)=A(c,d,S)$. One of such instances is
described in the following theorem.

\begin{thm}\label{thm-bcp}
If there exists a function $f(t,x)\in C^{1,2}([0,\infty)
\times\IR)$, such that $Y_t:= f(t,X_t)$ satisfies $dY_t =
\tilde\sigma_tdW_t$, where $\tilde\sigma_t\in C([0,\infty))$ is a
real, deterministic function satisfying $\tilde\sigma_t\neq 0$,
$\forall t\in[0,\infty)$, then there exists a standard BM $\{\tilde
W_s, s\geq 0\}$, such that for any boundaries $a(t), b(t)$,
\begin{enumerate}
\item
$P_X(a, b, T) = P_{\tilde W}(c, d, S)$, if $\tilde\sigma_t>0$,
$\forall t\in[0,\infty)$;
\item
$P_X(a, b, T) = P_{\tilde W}(d, c, S)$, if $\tilde\sigma_t<0$,
$\forall t\in[0,\infty)$;
\end{enumerate}
where
\begin{equation}\label{cs-general}
c(s) = f(t(s), a(t(s))) - f(0, x_0), 0\leq s\leq S,
\end{equation}
\begin{equation}\label{ds-general}
d(s) = f(t(s), b(t(s))) - f(0, x_0), 0\leq s\leq S,
\end{equation}
$t(s)$ is the inverse function of $s(t)=\int_0^t
\tilde\sigma_u^2du$ and $S=s(T)$.
\end{thm}

\proof%
First, from $dY_t =\tilde\sigma_tdW_t$ it follows that $\{Y_t,
t\geq 0\}$ is a continuous Gaussian process with variance function
\[
V(Y_t) = \int_0^t\tilde\sigma_u^2du = s(t), t\geq 0.
\]
Because $\tilde\sigma_t^2>0$, $s(t)$ is strictly increasing for
$t\geq 0$ and $s(0)=0$. For every $s\geq 0$, define $\tilde W_s:=
Y_{t(s)}-Y_0$, where $t(s)$ is the inverse function of $s(t)$.
Then since $\tilde W_0 = 0$ and
\[
V(\tilde W_s) = V(Y_{t(s)}) = s(t(s)) = s,
\]
$\{\tilde W_s, s\geq 0\}$ is a standard BM. Furthermore, for every
$t\geq 0$, $Y_t = Y_0 + \tilde W_{s(t)}$.

On the other hand, Assumption \ref{ass1} and It\^o's formula
(\ref{Itosform}) imply that
\[
\sigma(t,x)\frac{\partial f}{\partial x} = \tilde\sigma_t.
\]
Note that for every $t>0$,
\[
\frac{\partial f}{\partial x} = \frac{\tilde\sigma_t}{\sigma(t,x)}
\]
is continuous in $x$. It follows that, for any $x_0$ and $x$ in
the state space of $\{X_t\}$,
\begin{equation}\label{eqn-f}
f(t,x) = \tilde\sigma_t\int_{x_0}^x \frac{1}{\sigma(t,u)}du +
h(t),
\end{equation}
where $h(t)$ is an arbitrary function which does not depend on
$x$. Further, since $\tilde\sigma_t\neq 0$ and $\sigma(t,x)>0$,
for every $t>0$, $f(t,x)$ is a monotone function of $x$. For every
$t>0$, let $g(t,x)$ be the inverse function of $f(t,x)$. Then $X_t
= g(t, Y_t)$.

Now we consider case 1, where $\tilde\sigma_t>0$, $\forall
t\in[0,\infty)$. In this case function $f(t,x)$ is strictly
increasing in $x$. It follows that
\begin{eqnarray*}
P\left(a(t) < X_t < b(t), \forall t\in [0, T]\right) &=&
P\left(a(t) < g(t,Y_0 + \tilde W_{s(t)}) < b(t), \forall t\in [0, T]\right) \\
&=& P\left(f(t, a(t)) < Y_0 + \tilde W_{s(t)} < f(t, b(t)),
\forall t\in [0, T]\right) \\
&=& P\left(c(s) < \tilde W_s < d(s), \forall s\in [0, S]\right),
\end{eqnarray*}
where $c(s)$ and $d(s)$ are given by (\ref{cs-general}) and
(\ref{ds-general}) respectively, and $S=s(T)$.

Similarly, case 2 follows from the fact that function $f(t,x)$ is
strictly decreasing in $x$.
\endproof

\begin{rem}\label{rem-zeros}{\rm
The above derivation is readily generalized to more complicated case
where $\tilde\sigma_t$ has finite number of zeros. Without loss of
generality, suppose $\tilde\sigma_t$ has $k$ zeros
$0<t_1<t_2<\cdots<t_k<T$. Then $\tilde\sigma_t$ changes signs over
intervals $(t_{i-1}, t_i]$, $i=1,2,...,k+1$, where $t_0=0$ and
$t_{k+1}=T$. Denote by $T_+$ the union of all intervals in which
$\tilde\sigma_t>0$, and by $T_-$ the union of all intervals in which
$\tilde\sigma_t<0$. Then by (\ref{eqn-f}) function $f(t,x)$ is
strictly increasing for $t\in T_+$ and decreasing for $t\in T_-$.
Furthermore, let $S_+ = s(T_+)$, $S_-=s(T_-)$, where $s(t)=\int_0^t
\tilde\sigma_u^2du$. Then for any boundaries $a(t), b(t)$, we have
$P_X(a, b, T) = P_{\tilde W}(c, d, S)$, where
\[
c(s) = \left\{
\begin{array}{ll}
f(t(s), a(t(s))) - f(0,x_0), & \mbox{when $s\in S_+$} \\
f(t(s), b(t(s))) - f(0,x_0), & \mbox{when $s\in S_-$}
\end{array}
\right.
\]
\[
d(s) = \left\{
\begin{array}{ll}
f(t(s), b(t(s))) - f(0,x_0), & \mbox{when $s\in S_+$} \\
f(t(s), a(t(s))) - f(0,x_0), & \mbox{when $s\in S_-$}
\end{array}
\right.
\]
and $S=s(T)$.}
\end{rem}

There remains the question that for which processes $\{X_t\}$ do
functions $f(t,x)$ and $\tilde\sigma_t$ in Theorem \ref{thm-bcp}
exist. A class of such diffusion processes is characterized by the
next Theorem.

\begin{thm}\label{thm-pde}
Let $\{X_t, t\geq 0\}$ satisfy equation (\ref{sde-general}) with
drift $\mu(t,x)$ and diffusion coefficient $\sigma(t,x)\in
C^{1,2}([0,\infty)\times\IR)$. Then there exist functions $f(t,x)$
and $\tilde\sigma_t$ satisfying conditions of Theorem \ref{thm-bcp},
if and only if $\mu_t=\mu(t,x)$ and $\sigma_t=\sigma(t,x)$ satisfy
the following partial differential equation
\begin{equation}\label{pde-main}
\frac{\partial}{\partial x}\left[\frac{1}{\sigma_t}
\frac{\partial\sigma_t}{\partial t} + \sigma_t
\frac{\partial}{\partial x} \left(\frac{1}{2}
\frac{\partial\sigma_t}{\partial x} - \frac{\mu_t}{\sigma_t}
\right)\right]=0.
\end{equation}
\end{thm}
\proof%
Suppose $Y_t = f(t,X_t)$ and $dY_t = \tilde\sigma_tdW_t$, where
$f(t,x)$ and $\tilde\sigma_t$ satisfy conditions in Theorem
\ref{thm-bcp}. Then It\^o's formula (\ref{Itosform}) implies
\begin{equation}\label{Ito1}
\frac{\partial f}{\partial t} + \mu_t\frac{\partial f}{\partial x}
+ \frac{\sigma_t^2}{2} \frac{\partial^2f}{\partial x^2} = 0
\end{equation}
and
\begin{equation}\label{Ito2}
\sigma_t\frac{\partial f}{\partial x} = \tilde\sigma_t.
\end{equation}
From (\ref{Ito2}), we have
\[
\frac{\partial^2f}{\partial x^2} = -
\frac{\tilde\sigma_t}{\sigma_t^2}\frac{\partial\sigma_t}{\partial
x}.
\]
Substituting the above equation into (\ref{Ito1}), we obtain
\[
\frac{\partial f}{\partial t} = \tilde\sigma_t\left( \frac{1}{2}
\frac{\partial\sigma_t}{\partial x} -
\frac{\mu_t}{\sigma_t}\right).
\]
It follows that
\begin{equation}\label{f-def}
f(t,x) = \int_0^t\tilde\sigma_u\left(\frac{1}{2}
\frac{\partial\sigma_u}{\partial x} -
\frac{\mu_u}{\sigma_u}\right)du + \phi(x),
\end{equation}
where $\phi(x)\in C^2(\IR)$. Differentiating (\ref{f-def}) with
respect to $x$ and using equation (\ref{Ito2}), we have
\begin{equation}\label{phi-def}
\frac{\partial}{\partial x}\int_0^t\tilde\sigma_u\left(\frac{1}{2}
\frac{\partial\sigma_u}{\partial x} -
\frac{\mu_u}{\sigma_u}\right)du + \frac{d\phi(x)}{dx} =
\frac{\tilde\sigma_t}{\sigma_t}.
\end{equation}
Further, differentiating (\ref{phi-def}) with respect to $t$
yields
\[
\tilde\sigma_t\frac{\partial}{\partial x}\left(\frac{1}{2}
\frac{\partial\sigma_t}{\partial x} - \frac{\mu_t}{\sigma_t}
\right) = \frac{\partial}{\partial t}
\frac{\tilde\sigma_t}{\sigma_t},
\]
which is equivalent to
\begin{equation}\label{sigma-def}
\frac{1}{\sigma_t}\frac{\partial\sigma_t}{\partial t} +
\sigma_t\frac{\partial}{\partial x}\left(\frac{1}{2}
\frac{\partial\sigma_t}{\partial x} - \frac{\mu_t}{\sigma_t}
\right) = \frac{1}{\tilde\sigma_t}
\frac{\partial\tilde\sigma_t}{\partial t}.
\end{equation}
Since the right-hand side of (\ref{sigma-def}) does not depend on
$x$, differentiating both sides of this equation with respect to
$x$ yields (\ref{pde-main}).

Conversely, if (\ref{pde-main}) holds, then equation
(\ref{sigma-def}) defines a function $\tilde\sigma_t\in
C^1([0,\infty))$. Further, equation (\ref{phi-def}) defines a
function $\phi(x)\in C^2(\IR)$, whereas (\ref{f-def}) defines
function $f(t,x)$. From (\ref{f-def}) and (\ref{pde-main}), it is
easy to see that $f(t,x)$ satisfies condition in Theorem
\ref{thm-bcp}.
\endproof

Condition (\ref{pde-main}) is practical in real applications. Once
a process is given, one can easily check whether condition
(\ref{pde-main}) holds or not. On the other hand, equation
(\ref{pde-main}) can also be used to characterize certain
subclasses of diffusion processes, which may be of interest of
researchers working on particular problems. For example, for any
given type of diffusion coefficient $\sigma_t$, one can find all
possible functions $\mu_t$, such that (\ref{pde-main}) is
satisfied. Some examples are given in the following Corollary.

\begin{cor}\label{cor-subclasses}
For each of the following types of diffusion coefficient
$\sigma_t$, the corresponding general solutions of equation
(\ref{pde-main}) for the drift $\mu_t$ are given as follows, where
$\alpha(t)$, $\beta(t)$ are arbitrary real functions and
$\tilde\sigma(t): [0,\infty)\to\IR_+$.
\begin{enumerate}
\item
(L-class). For $\sigma(t,x) = \tilde\sigma(t)$, $\mu(t,x) =
\alpha(t)x + \beta(t)$, $x\in\IR$.
\item
(G-class). For $\sigma(t,x) = \tilde\sigma(t)x$, $\mu(t,x) =
\alpha(t)x + \beta(t)x\log x$, $x\in\IR_+$.
\end{enumerate}
\end{cor}

The proof of these results is straightforward and is therefore
omitted. It is easy to see that L-class corresponds to general
Ornstein-Uhlenbeck processes with time-dependent coefficients,
whereas G-class corresponds to growth processes which are widely
used in population genetics. It is easy to see that the geometric BM
process also belongs to G-class.

There remains a practical question as how to find transformations
$f(t,x)$ and $\tilde\sigma_t$ for a given process. Obviously this
depends on the functional forms of $\mu_t$ and $\sigma_t$ and can
only be worked out explicitly on case by case bases. Generally
speaking, equation (\ref{sigma-def}) defines $\tilde\sigma_t$, up
to a multiplicative constant, and (\ref{eqn-f}) defines $f(t,x)$
up a term $h(t)$ not depending on $x$, which is in turn determined
by (\ref{Ito1}). In the subsequent sections, we demonstrate this
for some widely used special processes and use the theoretical
results of this section to derive more detailed results for these
processes.

\begin{rem}\label{rem-transform}{\rm
The problem of obtaining the transition probability of a diffusion
process through one-to-one (strictly monotone) transformation of
Brownian motion was first posed by Kolmogorov (1931). A class of
such diffusion processes was characterized by Cherkasov (1957) and
Ricciardi (1976), whereas the most general class was derived by
Bluman (1980). Using this transformation method, Ricciardi,
Sacerdote and Sato (1984) obtained the FPT densities for a class of
diffusion processes crossing one-sided, continuously differentiable
boundaries. Since the main concern of this paper is the boundary
crossing probabilities, not the transformation itself, we do not
intend to be most general. The class of processes characterized by
Theorem \ref{thm-pde} is large enough to contain many interesting
processes arising in practice. Furthermore, the condition
(\ref{pde-main}) is practical and easy to work with in real
applications.}
\end{rem}

\section{Ornstein-Uhlenbeck processes}

Ornstein-Uhlenbeck (O-U) processes are a class of important
diffusion processes with wide applications. We start with the O-U
process which is defined in state space $\IR$ and satisfies
\begin{equation}\label{sde-ou}
dX_t =  \kappa(\alpha-X_t)dt + \sigma dW_t,\quad X_0 = x_0,
\end{equation}
where $\kappa,\sigma \in\IR_+$ and $\alpha\in\IR$ are constants. In
mathematical finance, (\ref{sde-ou}) is known as Vasicek model for
the short-term interest rate process (Vasicek 1977).

In this case, special forms of transformation $f(t,x)$ and
boundaries (\ref{cs-general}) and (\ref{ds-general}) can be derived.
For every $t>0$, if we define
\begin{equation}\label{eqn4}
Y_t = e^{\kappa t}(X_t-\alpha),\quad Y_0 = x_0-\alpha,
\end{equation}
then It\^o's formula and equation (\ref{sde-ou}) imply
\begin{eqnarray*}
dY_t &=& \kappa e^{\kappa t}(X_t-\alpha)dt + e^{\kappa t}dX_t \\
&=& \sigma e^{\kappa t}dW_t.
\end{eqnarray*}
It follows that $\{Y_t, t\geq 0\}$ is a continuous Gaussian process
with variance function
\begin{eqnarray*}
V(Y_t) &=& \sigma^2\int_0^te^{2\kappa u}du \\
&=& \frac{\sigma^2}{2\kappa}\left(e^{2\kappa t}-1\right) \\
&:=& s(t).
\end{eqnarray*}
Note that the function $s(t), t\geq 0$ is strictly increasing and
its inverse is given by
\[
t(s) = \frac{1}{2\kappa}\log\left(1+ \frac{2\kappa s}
{\sigma^2}\right), s\geq 0.
\]

For every $s\geq 0$, define $\tilde W_s:= Y_{t(s)}-Y_0$. Then, it
follows from $\tilde W_0 = 0$ and
\[
V(\tilde W_s) = V(Y_{t(s)}) = s(t(s)) = s,
\]
that $\{\tilde W_s, s\geq 0\}$ is a standard BM. Furthermore, for
every $t\geq 0$, $Y_t = Y_0 + \tilde W_{s(t)}$ and by (\ref{eqn4})
we have
\begin{eqnarray*}
X_t &=& e^{-\kappa t}Y_t + \alpha \nonumber\\
&=& e^{-\kappa t}\left(x_0-\alpha + \tilde W_{s(t)}\right) +\alpha.
\end{eqnarray*}
Therefore, we can write the one-sided BCP for $\{X_t\}$ as
\begin{eqnarray*}
P\left(X_t<b(t), \forall t\in [0, T]\right)%
&=& P\left(e^{-\kappa t}\left(x_0-\alpha + \tilde W_{s(t)}\right)
+ \alpha < b(t), \forall t\in [0, T]\right) \\
&=& P\left(\tilde W_{s(t)} < \alpha-x_0+e^{\kappa t}(b(t)-\alpha),
\forall t\in [0, T]\right) \\
&=& P\left(\tilde W_s<\alpha-x_0+e^{\kappa t(s)}(b(t(s))-\alpha),
\forall s\in[0, S]\right),
\end{eqnarray*}
where $S = s(T) = \sigma^2\left(e^{2\kappa T}-1\right)/2\kappa$.
Obviously, similar relation for the two-sided BCP can be established
in the same way. Thus, we have the following result.

\begin{cor}\label{cor-ou1}
Let $\{X_t, t\geq 0\}$ be an O-U process satisfying (\ref{sde-ou}).
Then there exists a standard BM $\{\tilde W_s, s\geq 0\}$, such that
for any boundaries $a(t), b(t)$,
\begin{equation}\label{bcp-ou}
P_X(a, b, T) = P_{\tilde W}(c, d, S),
\end{equation}
where
\begin{equation}\label{cs-ou}
c(s) = \alpha - x_0 + [a(t(s))-\alpha]\left(1+ \frac{2\kappa s}
{\sigma^2} \right)^{1/2},
\end{equation}
\begin{equation}\label{ds-ou}
d(s) = \alpha - x_0 + [b(t(s))-\alpha]\left(1+ \frac{2\kappa s}
{\sigma^2}\right)^{1/2},
\end{equation}
\begin{equation}\label{ts-ou}
t(s) = \frac{1}{2\kappa}\log\left(1+ \frac{2\kappa
s}{\sigma^2}\right), s\geq 0
\end{equation}
and $S = \sigma^2\left(e^{2\kappa T}-1\right)/2\kappa$.
\end{cor}

The above derivation can be extended to more general O-U processes
with time-dependent coefficients, i.e., processes which satisfy
\begin{eqnarray}\label{sde-ou2}
dX_t = \kappa(t)(\alpha(t) -X_t)dt+\sigma(t)dW_t,\quad X_0 = x_0,
\end{eqnarray}
where $\kappa(t), \sigma(t): [0,\infty)\mapsto\IR_+$ and $\alpha(t):
[0,\infty)\mapsto\IR$ are real, deterministic functions. Indeed, if
we define, for every $t>0$,
\[
Y_t = \exp\left(\int_0^t\kappa(u)du\right)(X_t-\gamma(t)),\quad Y_0
= x_0 - \gamma(0),
\]
where
\[
\gamma(t) = e^{-t}\alpha(0) + e^{-t}\int_0^te^u\alpha(u)du,
\]
then by It\^o's lemma, we have
\[
dY_t = \exp\left(\int_0^t\kappa(v)dv\right)\sigma(t)dW_t.
\]
It follows that
\begin{eqnarray}\label{eqn10}
s(t) = V(Y_t) =
\int_0^t\exp\left(2\int_0^u\kappa(v)dv\right)\sigma^2(u)du,
\end{eqnarray}
which is strictly increasing and hence admits the inverse function
$t(s), s\geq 0$. Therefore,
\begin{eqnarray*}
X_t &=& \gamma(t) + \exp\left(-\int_0^t\kappa(u)du\right)Y_t \\
&=& \gamma(t) + \exp\left(-\int_0^t\kappa(u)d\right)
\left(x_0-\gamma(0)+\tilde W_{s(t)}\right),
\end{eqnarray*}
where $\tilde W_{s} = Y_{t(s)} - Y_0$ is a standard BM. Thus, we
have the following result.
\begin{cor}\label{cor-ou2}
Let $\{X_t, t\geq 0\}$ be an O-U process satisfying (\ref{sde-ou2}).
Then there exists a standard BM $\{\tilde W_s, s\geq 0\}$, such that
for any boundaries $a(t), b(t)$, $P_X(a, b, T) = P_{\tilde W}(c, d,
S)$, where
\[
c(s) = \gamma(0)-x_0 + [a(t(s))-\gamma(t(s))]
\exp\left(\int_0^{t(s)}\kappa(u)du\right),
\]
\[
d(s) = \gamma(0)-x_0 + [b(t(s))-\gamma(t(s))]
\exp\left(\int_0^{t(s)}\kappa(u)du\right),
\]
$t(s)$ is the inverse function of $s(t)$ defined in (\ref{eqn10})
and $S=s(T)$.
\end{cor}

Relation (\ref{bcp-ou}) allows us to derive exact formulas for the
BCP of an O-U process $\{X_t\}$ crossing ceratin boundaries, using
the existing results for Brownian motion crossing, e.g., linear or
square-root boundaries. For the simplicity of notation, in the
following we present three examples for one-sided BCP.

\begin{exm}\label{exm-ou1}{\rm
First, consider boundaries $a(t)=-\infty$ and $b(t) = \alpha +
he^{\kappa t}, 0\leq t\leq T$, where $h\in\IR$ is an arbitrary
constant. Then by (\ref{cs-ou}) and (\ref{ds-ou}), $c(s)=-\infty$
and $d(s) = \alpha -x_0 + h + 2h\kappa s/\sigma^2$ is a linear
function in $s\in [0, S]$, where $S = \sigma^2(e^{2\kappa
T}-1)/2\kappa$. Therefore, by (\ref{bcp-ou}) and the well-known
formula for the BCP for BM crossing a linear boundary (e.g.,
equation (3) of Wang and P\"otzelberger 1997), we have
\[
P_X(a,b,T) = \Phi\left(\frac{he^{2\kappa T}+\alpha-x_0}
{\sigma\sqrt{(e^{2\kappa T}-1)/2\kappa}}\right) -
\exp\left(-\frac{4h\kappa(h+\alpha-x_0)}{\sigma^2}\right)
\Phi\left(\frac{he^{2\kappa T}-\alpha+x_0-2h}
{\sigma\sqrt{(e^{2\kappa T}-1)/2\kappa}}\right),
\]
where $\Phi$ is the standard normal distribution function.}
\end{exm}

\begin{exm}\label{exm-ou2}{\rm
Another set of boundaries are $a(t)=-\infty$ and $b(t) = \alpha +
he^{-\kappa t}, 0\leq t\leq T$. In this case $c(s)=-\infty$ and
$d(s) = \alpha -x_0 + h$ is a constant. Therefore, again by
(\ref{bcp-ou}) we have
\[
P_X(a,b,T) = 2\Phi\left(\frac{\alpha-x_0+h}{\sigma \sqrt{(e^{2\kappa
T}-1)/2\kappa}}\right) - 1.
\]
}
\end{exm}

\begin{exm}\label{exm-ou3}{\rm
Further, consider constant boundaries $a(t)=-\infty$ and $b(t) = h,
0\leq t\leq T$. Then $c(s)=-\infty$ and
\[
d(s) = \alpha -x_0 + (h-\alpha)\left(1+\frac{2\kappa
s}{\sigma^2}\right)^{1/2}
\]
is a square-root boundary, for which the "exact" BCP for BM is known
(Daniels 1996). A numerical example of this case is given in section
7.}
\end{exm}

\section{Growth processes}

Another important stochastic model in population genetics is the
growth processes (Ricciardi et al 1999), which is defined by the
SDE
\begin{equation}\label{sdf-growth}
dX_t = (\alpha X_t-\beta X_t\log X_t)dt+\sigma X_tdW_t,\quad X_0=
x_0
\end{equation}
and has state space $\IR_+$, where $\alpha, \beta$ and $\sigma$
are positive constants. It is easy to see that the growth
processes belong to the G-class described in Corollary
\ref{cor-subclasses}.

Similar to the O-U process, again by equations (\ref{Ito1}) and
(\ref{Ito2}), functions $f(t,x)$ and $\tilde\sigma_t$ can be
determined as
\[
f(t,x) = \frac{e^{\beta t}}{\sigma}\left(\log x +
\frac{\sigma^2-2\alpha}{2\beta}\right)
\]
and $\tilde\sigma_t = e^{\beta t}$. Hence the time transformation
is
\[
s(t) = \int_0^te^{2\beta u}du = \frac{1}{2\beta}(e^{2\beta t}-1),
t\geq 0.
\]
Therefore, by (\ref{cs-general}) and (\ref{ds-general}), we have
the following result.

\begin{cor}\label{cor-growth}
Let $\{X_t, t\geq 0\}$ be a growth process satisfying
(\ref{sdf-growth}). Then there exists a standard BM $\{\tilde W_s,
s\geq 0\}$, such that for any boundaries $a(t), b(t)$, $P_X(a, b,
T) = P_{\tilde W}(c, d, S)$, where
\begin{equation}\label{cs-growth}
c(s) = \frac{\sqrt{1+2\beta s}}{\sigma}\left(\log a(t(s)) +
\frac{\sigma^2-2\alpha}{2\beta}\right) - \frac{1}{\sigma}
\left(\log x_0 + \frac{\sigma^2-2\alpha}{2\beta}\right),
\end{equation}
\begin{equation}\label{ds-growth}
d(s) = \frac{\sqrt{1+2\beta s}}{\sigma}\left(\log b(t(s)) +
\frac{\sigma^2-2\alpha}{2\beta}\right) - \frac{1}{\sigma}
\left(\log x_0 + \frac{\sigma^2-2\alpha}{2\beta}\right),
\end{equation}
\begin{equation}\label{ts-growth}
t(s) = \frac{1}{2\beta}\log(1+2\beta s)
\end{equation}
and $S = (e^{2\beta T}-1)/2\beta$.
\end{cor}


Now we derive explicit formulas for the BCP for some special
boundaries.

\begin{exm}{\rm
Let $a(t)=0$ and
\[
b(t) = \exp\left(he^{\beta t} -
\frac{\sigma^2-2\alpha}{2\beta}\right), 0\leq t\leq T,
\]
where $h$ is a constant. Then by (\ref{cs-growth}) and
(\ref{ds-growth}), $c(s)=-\infty$ and
\[
d(s) = \frac{2h\beta s}{\sigma} + \frac{1}{\sigma} \left(h-\log
x_0 - \frac{\sigma^2-2\alpha}{2\beta}\right), 0\leq s\leq S
\]
which is a linear function and $S = (e^{2\beta T}-1)/2\beta$.
Therefore,
\begin{eqnarray*}
P_X(a,b,T) &=& \Phi\left(\frac{2\beta(he^{2\beta T}-\log x_0) -
\sigma^2+2\alpha}{\sigma\sqrt{2\beta(e^{2\beta T}-1)}}\right) \\
&-& \exp\left(\frac{4h\beta(\log x_0 -h)+2h(\sigma^2-2\alpha)}
{\sigma^2}\right) \Phi\left(\frac{2\beta(he^{2\beta T}- 2h + \log
x_0) + \sigma^2 - 2\alpha}{\sigma\sqrt{2\beta(e^{2\beta
T}-1)}}\right),
\end{eqnarray*}
where $\Phi$ is the standard normal distribution function.}
\end{exm}

\begin{exm}{\rm
Let $a(t)=0$ and
\[
b(t) = \exp\left(he^{-\beta t} -
\frac{\sigma^2-2\alpha}{2\beta}\right), 0\leq t\leq T,
\]
where $h$ is a constant. Then $c(s)=-\infty$ and
\[
d(s) = \frac{1}{\sigma} \left(h - \log x_0 -
\frac{\sigma^2-2\alpha}{2\beta}\right), 0\leq s\leq S.
\]
Therefore,
\[
P_X(a,b,T) = 2\Phi\left(\frac{2\beta(h-\log x_0) - \sigma^2 +
2\alpha}{\sigma\sqrt{2\beta(e^{2\beta T}-1)}}\right) - 1.
\]}
\end{exm}

\begin{exm}\label{exm-growth}{\rm
Again, consider the constant boundary $b(t) = h > 0, 0\leq t\leq T$.
Then
\[
d(s) = \frac{\sqrt{1+2\beta s}}{\sigma}\left(\log h +
\frac{\sigma^2-2\alpha}{2\beta}\right) - \frac{1}{\sigma}
\left(\log x_0 + \frac{\sigma^2-2\alpha}{2\beta}\right), 0\leq
s\leq S
\]
with $S = (e^{2\beta T}-1)/2\beta$. A numerical example of this
case is given in section 7.}
\end{exm}


\section{Geometric Brownian motion}

One popular stochastic model in mathematical finance is a
generalization of the classical Black-Scholes model (Black and
Scholes 1973) to the time-dependent interest rate process, as
defined by
\begin{equation}\label{sdf-gbm}
dX_t = r(t)X_tdt + \sigma X_tdW_t,\quad X_0=x_0,
\end{equation}
where $r(t): [0,\infty)\mapsto\IR_+$. It is well-known that in
this case there exists a risk-neutral probability measure, under
which the underlying process is a geometric BM
\[
X_t = X_0\exp\left(R(t) - \frac{\sigma^2t}{2} +\sigma\tilde W_t
\right),
\]
where $R(t) = \int_0^tr(u)du$ and $\{\tilde W_t\}$ is a standard BM
under the same risk-neutral probability measure. In the following we
use the approach of this paper to derive the BCP for $\{X_t\}$.

First, since equation (\ref{Ito2}) now becomes
\[
\frac{\partial f}{\partial x} = \frac{\tilde\sigma_t}{\sigma x},
\]
function $f(t,x)$ is given by
\[
f(t,x) = \frac{\tilde\sigma_t}{\sigma}\left(\log x + h(t)\right),
\]
where $h(t)$ is an arbitrary function to be determined later.
Further, since
\[
\frac{\partial f}{\partial t} = \frac{(\log x + h(t))}{\sigma}
\frac{\partial\tilde\sigma_t}{\partial t} +
\frac{\tilde\sigma_t}{\sigma}\frac{\partial h(t)}{\partial t},
\]
and
\[
\frac{\partial^2f}{\partial x^2} =
-\frac{\tilde\sigma_t}{\sigma x^2},
\]
equation (\ref{Ito1}) becomes
\[
(\log x + h(t))\frac{\partial\tilde\sigma_t}{\partial t} +
\tilde\sigma_t \left(\frac{\sigma^2}{2}- r(t)- \frac{\partial
h(t)}{\partial t}\right) = 0,
\]
which holds, if $\tilde\sigma_t=1$ and
\[
\frac{\partial h(t)}{\partial t} = \frac{\sigma^2}{2}- r(t).
\]
The last equation above has a particular solution
\begin{eqnarray*}
h(t) &=& \int_0^t\left(\frac{\sigma^2}{2}- r(u)\right)du \\
&=& \frac{\sigma^2t}{2} - R(t).
\end{eqnarray*}
Therefore, we obtain $\tilde\sigma_t=1$ and
\[
f(t,x) = \frac{1}{\sigma}\left(\log x + \frac{\sigma^2t}{2} - R(t)
\right).
\]
By Theorem \ref{thm-bcp} the time transformation is $s(t) = t, t\geq
0$, and by (\ref{cs-general}) and (\ref{ds-general}), we have the
following well-known result.

\begin{cor}\label{cor-gbm}
Let $\{X_t, t\geq 0\}$ be a geometric BM satisfying
(\ref{sdf-gbm}). Then there exists a standard BM $\{\tilde W_s,
s\geq 0\}$, such that for any boundaries $a(t), b(t)$, $P_X(a, b,
T) = P_{\tilde W}(c, d, T)$, where
\begin{equation}\label{cs-gbm}
c(t) = \frac{1}{\sigma}\left(\log\left(\frac{a(t)}{x_0}\right) +
\frac{\sigma^2t}{2} - R(t)\right), 0\leq t\leq T,
\end{equation}
and
\begin{equation}\label{ds-gbm}
d(t) = \frac{1}{\sigma}\left(\log\left(\frac{b(t)}{x_0}\right) +
\frac{\sigma^2t}{2} - R(t)\right), 0\leq t\leq T.
\end{equation}
\end{cor}


Now we consider some special cases.

\begin{exm}{\rm
Let $a(t) = 0$ and $b(t)= \exp(pt+q+R(t))$, where $p,q$ are
constants. Then by (\ref{cs-gbm}) and (\ref{ds-gbm}), $c(t) =
-\infty$ and
\[
d(t) = \frac{1}{\sigma}\left(p+\frac{\sigma^2}{2}\right)t +
\frac{q-\log x_0}{\sigma}, 0\leq t\leq T.
\]
Therefore,
\begin{eqnarray*}
P_X(a,b,T)
&=& \Phi\left(\frac{(p+\sigma^2/2)T+q-\log x_0}{\sigma\sqrt{T}}\right)\\
&-& \exp\left(\frac{(2p+\sigma^2)(\log x_0-q)}{\sigma^2}\right)
\Phi\left(\frac{(p+\sigma^2/2)T-q+\log x_0}{\sigma\sqrt{T}}
\right).
\end{eqnarray*}
}
\end{exm}

\begin{exm}{\rm
Another special case is the constant interest rate $r(t)=r$, so
that $R(t) = rt$. In this case, if we take constant boundaries
$a(t)=0$ and $b(t) = h>0$, then $c(t)=-\infty$ and
\[
d(t) = \frac{1}{\sigma}\left(\frac{\sigma^2}{2}-r\right)t +
\frac{1}{\sigma}\log\left(\frac{h}{x_0}\right), 0\leq t\leq T.
\]
Therefore,
\begin{eqnarray*}
P_X(a,b,T)
&=& \Phi\left(\frac{(\sigma^2/2-r)T+\log(h/x_0)}{\sigma\sqrt{T}}\right)\\
&-& \exp\left(\frac{(2r-\sigma^2)\log(h/x_0)}{\sigma^2}\right)
\Phi\left(\frac{(\sigma^2/2-r)T-\log(h/x_0)}{\sigma\sqrt{T}}
\right).
\end{eqnarray*}
}
\end{exm}

More numerical examples for the geometric and standard BM are
given in the next section.


\section{Numerical examples}

In this section we compute some numerical examples for the BCP for
various processes discussed in previous sections. For the simplicity
of notation, we present examples of one-sided BCP only.

Given any diffusion process $\{X_t\}$ and a boundary, the BCP for
$\{X_t\}$ can be given by the BCP for Brownian motion and a
transformed boundary $b$, i.e., $P_W(-\infty,b,T)$. As discussed
in section 2, the BCP $P_W(-\infty,b,T)$ can be approximated by
$P_W(-\infty,b_n,T)$, where $b_n$ is a piecewise linear boundary
converging to $b$ uniformly on the interval $[0,T]$. In
particular, if $b_n(t)$ and $\tilde b_n(t)$ are such piecewise
linear functions approaching $b(t)$ from below and above
respectively, then
\[
P_W(-\infty,b_n,T)\leq P_W(-\infty,b,T)\leq P_W(-\infty,\tilde
b_n,T).
\]
The two sides of the above inequalities are given by
(\ref{bcp-bm}) with function $g$ given in (\ref{bcp-one-sided}).
The expectation in (\ref{bcp-bm}) can be easily computed through
Monte Carlo simulation, because $W_{t_1}$, $W_{t_2}$, ...,
$W_{t_n}$ have a multivariate normal distribution.

For all examples in this section, $b_n(t)$ and $\tilde b_n(t)$ are
constructed with $n=128$ equally-spaced nodes and lower and upper
bounds for the corresponding BCP are computed. In each Monte Carlo
simulation, $10^6$ repetitions are carried out.

\begin{exm}[O-U process]{\rm
Consider constant boundaries of Example \ref{exm-ou3}, where
$a(t)=-\infty$ and $b(t) = h, 0\leq t\leq T$. Then, $c(s)=-\infty$
and $d(s)$ is a square-root boundary. In particular, we set
parameters to $\alpha=x_0$, $\sigma^2=2\kappa=1$, $h=x_0+1$ and
$T=1$. Then $d(s)=\sqrt{1+s}$, for $0\leq s\leq S=e-1$. Using the
procedure described above, the computed lower and upper bounds for
the BCP are $0.721463\leq P_X(a,b,T)\leq 0.721464$, with a
simulation standard error of $0.000440$.}
\end{exm}

\begin{exm}[Growth process]{\rm
Now consider boundaries of Example \ref{exm-growth}, where $b(t) =
h, 0\leq t\leq T$. Then $d(s)$ is again a square-root boundary. In
particular, if we take $\sigma^2=2\alpha, \beta=1/2$, $x_0=1$,
$h=e^\sigma$ and $T=1$, then $d(s)=\sqrt{1+s}$ and $S=e-1$, which is
the same as the square-root boundary of Example 20. The lower and
upper bounds for the BCP are therefore given there.}
\end{exm}

\begin{exm}[Geometric BM]{\rm
Roberts and Shortland (1997) give an example of pricing a European
call option, who's underlying security $X_t$ follows the SDE
(\ref{sdf-gbm}) with $\sigma = 0.1$, $r(t) = 0.1+0.05e^{-t}$ and
$x_0 = 10$. The option has a knock-in boundary $b(t)=12$ and
maturity $T=1$. Then by Corollary \ref{cor-gbm} and equation
(\ref{ds-gbm}), we have $P_X(-\infty,b,T) = P_{\tilde
W}(-\infty,d,T)$, where
\[
d(t) = 10\log(1.2)-0.5-0.95t+0.5e^{-t}, 0\leq t\leq T.
\]
Therefore, using the approach of this paper the lower and upper
bounds for the BCP are found to be $0.603728\leq
P_X(-\infty,b,T)\leq 0.603729$, with a simulation standard error
$0.000483$.}
\end{exm}

\begin{exm}[BM and Daniels boundary]{\rm
Finally, consider the well-known Daniels' boundary for standard BM
$W=\{W_t,t\geq 0\}$ (Daniels 1969, 1996), which is
\[
b(t)= \frac{1}{2} - t\log\left(\frac{1}{4}+\frac{1}{4}
\sqrt{1+8e^{-1/t}}\right), 0\leq t\leq T,
\]
with $T=1$. This boundary has been used by many authors to test
their computational algorithms. The "exact" BCP for this boundary is
known to be $P_W(-\infty,b,T)= 0.520251$. The procedure described
above gives an approximation for the corresponding BCP as
$P_W(-\infty,b,T)\approx 0.520293$ with a simulation standard error
$0.000490$.}
\end{exm}

\section{Conclusions and Discussion}

We proposed a direct method for computing boundary crossing
probabilities (BCP) for a class of diffusion processes which can be
expressed as piecewise monotone functionals of a standard Brownian
motion. This class includes many interesting processes in real
applications, e.g., Ornstein-Uhlenbeck, growth processes and
geometric Brownian motion with time dependent drift. This method
applies to both one-sided and two-sided general nonlinear
boundaries, which may be discontinuous. Using this approach explicit
formulas for boundary crossing probability for certain nonlinear
boundaries have been obtained, which are useful to evaluate and
compare computational algorithms.

Compared to the traditional approach for the first-passage time
densities, the BCP method can be applied to both one-sided and
two-sided boundaries. It also applies to discontinuous boundaries.
The actual computation can be easily done using Monte Carlo
simulation methods, which is practical and easy to implement.
Furthermore, the approximation errors are also computed
automatically. The approach of this paper is based on the stochastic
differential equation for the processes. It is different from the
traditional methods of first-passage time densities, which is based
on the partial differential equation for the transition density
function.


\section*{Acknowledgment}

We would like to thank the editor Joseph Glaz and two anonymous
referees for their helpful comments and suggestions. Financial
support from the Natural Sciences and Engineering Research Council
of Canada (NSERC) is gratefully acknowledged.


\end{document}